\documentstyle[twoside,amstex]{article}
\input amssym.def
\input amssym.tex
\def\bc{\begin{center}}
\def\ec{\end{center}}
\def\no{\noindent}
\def\hang{\hangindent\parindent}
\def\textindent#1{\indent\llap{[#1]\enspace}\ignorespaces}
\def\re{\par\hang\textindent}
\begin{document}
\thispagestyle{empty} \vspace*{3 true cm} \pagestyle{myheadings}
\markboth {\hfill {\sl Huanyin Chen}\hfill} {\hfill{\sl ON
UNIQUELY $\pi$-CLEAN RINGS}\hfill} \vspace*{-1.5 true cm}
\bc{\large\bf ON UNIQUELY $\pi$-CLEAN RINGS}\ec

\vskip6mm
\bc{{\bf Huanyin Chen}\\[1mm]
Department of Mathematics, Hangzhou Normal University\\
Hangzhou 310036, People's Republic of China}\ec

\vskip10mm
\begin{minipage}{120mm}
\no {\bf Abstract:} An element of a ring is unique clean if it can
be uniquely written as the sum of an idempotent and a unit. A ring
$R$ is uniquely $\pi$-clean if some power of every element in $R$
is uniquely clean. In this article, we prove that a ring $R$ is
uniquely $\pi$-clean if and only if for any $a\in R$, there exists
an $m\in {\Bbb N}$ and a central idempotent $e\in R$ such that
$a^m-e\in J(R)$, if and only if $R$ is abelian; every idempotent
lifts modulo $J(R)$; and $R/P$ is torsion for all prime ideals $P$
containing the Jacobson radical $J(R)$. Further, we prove that a
ring $R$ is uniquely $\pi$-clean and $J(R)$ is nil if and only if
$R$ is an abelian periodic ring, if and only if for any $a\in R$,
there exists some $m\in {\Bbb N}$ and a unique idempotent $e\in R$
such that $a^m-e\in P(R)$, where $P(R)$ is the prime radical of
$R$. \vskip3mm \no {\bf MR(2010) Subject Classification}: 16U99,
16E50.
\end{minipage}

\vskip15mm \bc{\bf 1. INTRODUCTION}\ec

\vskip4mm \no An element of a ring is unique clean if it can be
uniquely written as the sum of an idempotent and a unit. A ring
$R$ is uniquely clean if every element in $R$ is uniquely clean.
Many results on such rings can be found in [2], [4] and [10].
Following Zhou [10], a ring $R$ is uniquely $\pi$-clean if some
power of every element in $R$ is uniquely clean. This is a natural
generalization of uniquely clean rings. Some structures of such
rings was claimed in [10]. The motivation of this paper is to
develop further characterizations of uniquely $\pi$-clean rings.

In Section 2, we investigate the structure theorems of uniquely
$\pi$-clean rings, and prove that a ring $R$ is uniquely
$\pi$-clean if and only if for any $a\in R$, there exists an $m\in
{\Bbb N}$ and a central idempotent $e\in R$ such that $a^m-e\in
J(R)$, if and only if for any $a\in R$, there exists an $m\in
{\Bbb N}$ and a unique $e\in R$ such that $a^m-e\in J(R)$, and
$J(R)=\{ x\in R~|~x^m-1\in U(R)~\mbox{for all}~m\in {\Bbb N}\}.$
An exchange-like characterization of such rings are also obtained.
In Section 3, we characterize uniquely $\pi$-cleanness by means of
some prime ideals. It is shown that a ring $R$ is strongly
$\pi$-clean if and only if $R$ is abelian; every idempotent lifts
modulo $J(R)$; and $R/P$ is torsion for all prime ideals $P$
containing the Jacobson radical $J(R)$. Furthermore, we consider a
type of radical-like ideal $J^*(R)$, and characterize uniquely
$\pi$-clean ring $R$ by using such special one. A ring $R$ is
periodic if for any $a\in R$ there exist distinct $m,n\in {\Bbb
N}$ such that $a^m=a^n$. In the last section, we investigate
uniquely $\pi$-clean ring with nil Jacobson radical. We prove that
a ring $R$ is uniquely $\pi$-clean and $J(R)$ is nil if and only
if $R$ is an abelian periodic ring, if and only if for any $a\in
R$ there exists some $m\in {\Bbb N}$ such that $a^m\in R$ is
uniquely nil clean, if and only if for any $a\in R$, there exists
some $m\in {\Bbb N}$ and a unique idempotent $e\in R$ such that
$a^m-e\in P(R)$, where $P(R)$ is the prime radical of $R$. Here,
an element $a\in R$ is uniquely nil clean if there exists a unique
idempotent $e\in R$ such that $a-e\in R$ is nilpotent ([2] and
[4]).

Throughout, all rings are associative with an identity. We use
$J(R)$ and $P(R)$ to denote the the Jacobson radical and prime
radical of a ring $R$. $N(R)$ stands for the set of all nilpotent
elements in $R$.

\vskip10mm\bc{\bf 2. STRUCTURE THEOREMS}\ec

\vskip4mm The aim of this is to explore the structures of uniquely
$\pi$-clean rings. Recall that a ring $R$ is an exchange ring if
for any $a\in R$ there exists an idempotent $e\in aR$ such that
$1-e\in (1-a)R$. A ring $R$ is an exchange ring if and only if,
for every right $R$-module $A$ and any two decompositions
$A=M\oplus N=\bigoplus_{i\in I} A_i$, where $M_R\cong R$ and the
index set $I$ is finite, there exist submodules $A_i'\subseteq
A_i$ such that $A=M\oplus \big(\bigoplus_{i\in I} A_i'\big)$. The
class of exchange rings is very large. For instances, regular
rings, $\pi$-regular rings, strongly $\pi$-regular rings,
semiperfect rings, left or right continuous rings, clean rings and
unit $C^*$-algebras of real rank zero, etc. We begin with

\vskip4mm \hspace{-1.8em} {\bf Lemma 2.1.}\ \ {\it Every uniquely
$\pi$-clean ring is an abelian exchange
ring.}\vskip2mm\hspace{-1.8em} {\it Proof.}\ \ Let $R$ be uniquely
$\pi$-clean, let $e\in R$ be an idempotent, and let $r\in R$.
Choose $x=1-\big(e+er(1-e)\big)$. Then there exists some $n\in
{\Bbb N}$ such that $x^n\in R$ is uniquely clean. One easily
checks that
$x^n=x=e+\big(1-2e-er(1-e)\big)=\big(e+er(1-e)\big)+\big(1-2(e+er(1-e))\big)$.
Further, $e=e^2\in R,
\big(1-2e+er(1-e)\big)^{-1}=\big(1-er(1-e)\big)(1-2e),
\big(e+er(1-e)\big)=\big(e+er(1-e)\big)^2$ and
$\big(1-2(e+er(1-e))\big)^2=1$. By the uniqueness, we get
$e=e+er(1-e)$, and then $er=ere$. Likewise, $re=ere$. Thus,
$er=re$, and therefore $R$ is abelian.

For any $a\in R$, then we can find some $m\in {\Bbb N}$ such that
$a^m\in R$ is clean. Write $a^m=f+v$, where $f=f^2,v\in U(R)$.
Then $a^m-f^m=v$, and so $a-f\in U(R)$. This implies that $R$ is
strongly clean. In view of [9, Theorem 30.2], $R$ is an exchange
ring.\hfill$\Box$

\vskip4mm A ring $R$ is strongly clean if for any $a\in R$ there
exists an idempotent $e\in R$ such that $a-e\in U(R)$ and $ea=ae$.
As a consequence of Lemma 2.1, every uniquely $\pi$-clean ring is
strongly clean. A ring $R$ is uniquely clean provided that every
element in $R$ can be uniquely written as the sum of an idempotent
and a unit. It is easy to verify that ${\Bbb Z}/3{\Bbb Z}$ is not
uniquely clean as $2=0+2=1+1$, while ${\Bbb Z}/3{\Bbb Z}$ is
uniquely $\pi$-clean. Let $R=\bigoplus\limits_{p~\mbox{is
prime}}{\Bbb Z}/(p+1){\Bbb Z}$. Then $R$ is strongly clean. For
any $1\leq m\leq \big[ log_2p\big], 2^m\in {\Bbb Z}/(p+1){\Bbb Z}$
is not uniquely clean. Thus, $R$ is not uniquely $\pi$-clean.
Therefore, we conclude that $\{~\mbox{uniquely clean
rings}\}~\subsetneq \{~\mbox{uniquely}~ \pi\mbox{--clean
rings}~\}\subsetneq \{~\mbox{strongly clean rings}~\}.$

\vskip4mm \hspace{-1.8em} {\bf Theorem 2.2.}\ \ {\it Let $R$ be a
ring. Then $R$ is uniquely $\pi$-clean if and only if}
\begin{enumerate}
\item [(1)] {\it $R$ is abelian;}
\vspace{-.5mm}
\item [(2)] {\it Every idempotent lifts modulo $J(R)$;}
\vspace{-.5mm}
\item [(3)] {\it $R/J(R)$ is uniquely $\pi$-clean.}
\end{enumerate}
\vspace{-.5mm} {\it Proof.}\ \ Suppose $R$ is uniquely
$\pi$-clean. In view of Lemma 2.1, $R$ is an abelian exchange
ring. This proves $(1)$ and $(2)$, in terms of [9, Theorem 30.2].
For any $\overline{a}\in R/J(R)$, then $a\in R$ is uniquely
$\pi$-clean. Thus, we have some $n\in {\Bbb N}$ such that $a^n\in
R$ is uniquely clean. This implies that $a^n=e+u,e=e^2\in R, u\in
U(R)$. Hence, $\overline{a}^n=\overline{e}+\overline{u}$. Write
$\overline{a}^n=\overline{f}+\overline{v},
\overline{f}=\overline{f}^2\in R/J(R), \overline{v}\in
U\big(R/J(R)\big)$. Clearly, every unit lifts modulo $J(R)$. So we
may assume that $f=f^2\in R,v\in U(R)$. As a result, there exists
some $r\in J(R)$ such that $a^n=e+u=f+(v+r)$. By the uniqueness,
we get $e=f$. Therefore $R/J(R)$ is uniquely $\pi$-clean.

Conversely, assume that $(1)-(3)$ hold. For any $a\in R$, we have
$\overline{a}\in R/J(R)$, and so there exists some $n\in {\Bbb N}$
such that $\overline{a}^n\in R$ is uniquely clean. By hypothesis,
idempotents lift modulo $J(R)$. In addition, units lift modulo
$J(R)$. Thus, $a^n=e+u,e=e^2\in R, u\in U(R)$. Write $a^n=f+v,
f=f^2,v\in U(R)$. Then $\overline{a}^n=\overline{f}+\overline{v}$.
By the uniqueness, we get $\overline{e}=\overline{f}$, i.e.,
$e-f\in J(R)$. This infers that $f(1-e)=(e-f)(e-1)\in J(R)$. As
every idempotent in $R$ is central, $f(1-e)\in R$ is an
idempotent, thus, $f(1-e)=0$. It follows that $f=fe$. Likewise,
$e=ef$. Consequently, $e=f$, and therefore $R$ is uniquely
$\pi$-clean.\hfill$\Box$

\vskip4mm \hspace{-1.8em} {\bf Corollary 2.3.}\ \ {\it Every
corner of a uniquely $\pi$-clean ring is uniquely $\pi$-clean.}
\vskip2mm\hspace{-1.8em} {\it Proof.}\ \ Let $R$ be uniquely
$\pi$-clean, and let $e=e^2\in R$. In light of Theorem 2.2, $e\in
R$ central. For any $eae\in eRe$, then $eae+1-e\in R$ is uniquely
$\pi$-clean. So we have some $n\in {\Bbb N}$ such that
$(eae+1-e)^n\in R$ is uniquely clean. Thus, $(eae+1-e)^n=f+u,
f=f^2\in R, u\in U(R)$, and so $(eae)^n=efe+eue$ is clean in
$eRe$. Write $(eae)^n=g+v, g=g^2\in eRe, v\in U(eRe)$. Then
$(eae+1-e)^n=(eae)^n+1-e=g+(v+1-e)$, where $g=g^2\in R, v+1-e\in
U(R)$. Thus, $g=f=ege=efe$, as required.\hfill$\Box$

\vskip4mm Lemma 2.1 shows that every uniquely $\pi$-clean ring is
an abelian exchange ring. We next exhibit an exchange-like
property of such rings.

\vskip4mm \hspace{-1.8em} {\bf Theorem 2.4.}\ \ {\it Let $R$ be a
ring. Then $R$ is uniquely $\pi$-clean if and only
if\/}\vspace{-.5mm}
\begin{enumerate}
\item [(1)] {\it $R$ is abelian;}
\vspace{-.5mm}
\item [(2)] {\it For any $a\in R$, there exists an $n\in {\Bbb N}$ and a unique idempotent $e\in a^nR$ such that $1-e\in (1-a^{n})R$.}
\end{enumerate}
\vspace{-.5mm} {\it Proof.}\ \ Suppose that $R$ is uniquely
$\pi$-clean. In view of Lemma 2.1, every idempotent in $R$ is
central. For any $a\in R$, there exists some $n\in {\Bbb N}$ such
that $a^n\in R$ is uniquely clean. Write $a^n=f+v$, where
$f=f^2,v\in U(R)$. Set $g=v(1-f)v^{-1}$. Then $g=g^2\in R$.
Obviously, we get
$$\begin{array}{lll}
(a^n-g)v&=&\big(f+v-v(1-f)v^{-1}\big)v\\
&=&v^2+fv-v+vf\\
&=&a^{2n}-a^n. \end{array}$$ Thus $g-a^n\in (a^n-a^{2n})R$, and so
$g\in a^nR$ and $1-g\in (1-a^n)R$.

If there exists an idempotent $h\in a^nR$ such that $h\in a^nR$
and $1-h\in (1-a^n)R$. Write $h=a^nx, xh=x$. Then $xa^nx=x$. It is
easy to verify that $xa^n=x(a^nx)a^n=a^nx(xa^n)=a^n(xa^n)x=a^nx$.
Write $1-h=(1-a^n)y, y(1-h)=y$. Likewise, $y(1-a^n)=(1-a^n)y$. One
directly checks that $\big(a^n-(1-h)\big)^{-1}=x-y$, i.e.,
$a^n-(1-h)\in U(R)$, By the uniqueness, we get $1-h=f$. Hence,
$g=v(1-f)v^{-1}=1-f=h$, as desired.

Conversely, assume that $(1)$ and $(2)$ hold. For any $a\in R$,
there exists an $n\in {\Bbb N}$ and a unique idempotent $e\in
a^nR$ such that $1-e\in (1-a^{n})R$. As in the preceding
discussion, we get $a^n-(1-e)\in U(R)$. Write $a^n=f+v$, where
$f=f^2,v\in U(R)$. Set $g=v(1-f)v^{-1}$. Then $g=g^2\in R$.
Further, we have $g\in a^nR$ and $1-g\in (1-a^n)R$. By the
uniqueness, we obtain $g=e$, and so $v(1-f)v^{-1}=e$. Thus,
$f=1-e$, hence the result.\hfill$\Box$

\vskip4mm \hspace{-1.8em} {\bf Corollary 2.5.}\ \ {\it Let $R$ be
a ring. Then $R$ is uniquely $\pi$-clean if and only
if\/}\vspace{-.5mm}
\begin{enumerate}
\item [(1)] {\it Every idempotent in $R$ is central.}
\vspace{-.5mm}
\item [(2)] {\it For any $a\in R$, there exists an $n\in {\Bbb N}$ and a unique idempotent $e\in Ra^n$ such that $1-e\in R(1-a^{n})$.}
\end{enumerate}
\vspace{-.5mm} {\it Proof.}\ \ Obviously, a ring $R$ is uniquely
$\pi$-clean if and only if so is the opposite ring $R^{op}$.
Applying Theorem 2.4 to $R^{op}$, we complete the
proof.\hfill$\Box$

\vskip4mm A ring $R$ is local if it has only one maximal right
ideal. A ring $R$ is potent if for any $a\in R$ there exists some
$n\in {\Bbb N}$ such that $a^n=a$. We note that every potent ring
is commutative.

\vskip4mm \vskip4mm \hspace{-1.8em} {\bf Lemma 2.6.}\ \ {\it Let
$R$ be a local ring. If $R$ is uniquely $\pi$-clean, then $R/J(R)$
is potent.}\vskip2mm\hspace{-1.8em} {\it Proof.}\ \ Suppose that
there exists some $a\in R$ such that $a^n-a\not\in J(R)$ for all
$n\geq 2$. Then $a(a^{n-1}-1)\in U(R)$ as $R$ is a local ring.
This implies that $a\in U(R)$ and $a^{n-1}-1\in U(R)$ for all
$n\geq 2$. Since $R$ is uniquely $\pi$-clean, we have an $m\in
{\Bbb N}$ such that $a^m\in R$ is uniquely clean. But
$a^m=0+a^m=1+(a^m-1)$, a contradiction. Therefore, for any $a\in
R$, there exists some integer $n\geq 2$ such that $a^n-a\in J(R)$.
That is, $R/J(R)$ is potent.\hfill$\Box$

\vskip4mm The following lemma was firstly claimed by Lee and Zhou
[11] without proof. Here, we include an alternative proof for the
self-contained.

\vskip4mm \hspace{-1.8em} {\bf Lemma 2.7.}\ \ {\it Let $R$ be a
ring. Then $R$ is uniquely $\pi$-clean if and only if}
\begin{enumerate}
\item [(1)] {\it $R$ is abelian;}
\vspace{-.5mm}
\item [(2)] {\it Every idempotent lifts modulo $J(R)$;}
\vspace{-.5mm}
\item [(3)] {\it $R/J(R)$ is potent.}
\end{enumerate}
\vspace{-.5mm} {\it Proof.}\ \ Suppose that $R$ is uniquely
$\pi$-clean. In view of Theorem 2.2, $R/J(R)$ is uniquely
$\pi$-clean, idempotents lift modulo $J(R)$ and idempotents in $R$
are central. Clearly, $R/J(R)$ is isomorphic to a direct product
of some primitive rings $R_i$. Thus, $R_i$ is a homomorphic image
of $R/J(R)$. In view of Lemma 2.1, $R_i$ be an abelian exchange
ring, and that $R_i$ is strongly $\pi$-clean by an argument in
[11, Examples]. As every abelian exchange primitive ring is local,
$R_i$ is local. Therefore, $R_i$ is potent from Lemma 2.6. This
shows that $R/J(R)$ is potent, as desired.

Conversely, assume that $(1)-(3)$ hold. By hypothesis. $S:=R/J(R)$
is potent. Let $a\in S$. Then $a^m=a$ for some $m\geq 2$. Thus,
$a^{m-1}\in S$ is an idempotent. Hence,
$a^{m-1}=\big(1-a^{m-1}\big)+\big(2a^{m-1}-1)$, where
$1-a^{m-1}\in S$ is an idempotent and
$2a^{m-1}-1=\big(2a^{m-1}-1)^{-1}\in U(S)$. If there exist an
idempotent $f\in S$ and a unit $u\in R$ such that $a^{m-1}=f+u$.
As $R$ is abelian, it is easy to verify that
$\big(a^{m-1}+f-1\big)\big(a^{m-1}-f)^2=0$, and so $f=1-a^{m-1}$;
hence, $R/J(R)$ is strongly $\pi$-clean. Therefore we complete the
proof by Theorem 2.2.\hfill$\Box$

\vskip4mm \hspace{-1.8em} {\bf Theorem 2.8.}\ \ {\it Let $R$ be a
ring. Then the following are equivalent:}
\begin{enumerate}
\item [(1)] {\it $R$ is uniquely $\pi$-clean.}
\vspace{-.5mm}
\item [(2)] {\it For any $a\in R$, there exists an $m\in {\Bbb N}$ and a central idempotent $e\in R$ such that
$a^m-e\in J(R)$.}
\end{enumerate}
\vspace{-.5mm} {\it Proof.}\ \ $(1)\Rightarrow (2)$ In view of
Lemma 2.7, $R/J(R)$ is potent. For any $a\in R$, $\overline{a}\in
R/J(R)$ is potent, and so $\overline{a}^m\in R/J(R)$ is an
idempotent for some $m\in {\Bbb N}$. By using Lemma 2.7 again, we
can find a central idempotent $e\in R$ such that
$\overline{a}^m=\overline{e}$, and so $a^m-e\in J(R)$.

$(2)\Rightarrow (1)$ If $e\in R$ is an idempotent, then we have a
central idempotent $f\in R$ such that $e-f\in J(R)$. As
$(e-f)^3=e-f$, we deduce that $e=f$; hence, every idempotent in
$R$ is central. If $e-e^2\in J(R)$, then we can find a central
idempotent $f\in R$ such that $e^m-f\in J(R)$ for some $m\in {\Bbb
N}$. As $e-e^2\in J(R)$, if $m\geq 3$, we see that
$e-e^m=(e-e^2)+(e-e^2)e+\cdots +(e-e^2)e^{m-2}\in J(R)$. Thus
$e-f\in J(R)$, and then idempotents lift modulo $J(R)$.

For any $a\in R$, there exist $m\in {\Bbb N}$ such that $a^m-e\in
J(R)$ for a central idempotent. Hence,
$\overline{a}^m=\overline{e}$ in $R/J(R)$. Thus, $S:=R/J(R)$ is
periodic. Thus, $S$ is an abelian exchange ring. If $x^2=0$ and
$x\neq 0$ in $S$, then $x\not\in J(S)$. For any $r\in R$, there
exists some $g\in Rxr$ such that $1-g\in R(1-xr)$. Write $g=crx$
for a $c\in S$. Then $g=crgx=(cr)^2x^2=0$, and so $1-xr\in R$ is
left invertible. As $R$ is abelian, it is easy to check that
$1-xr\in U(R)$. This shows that $x\in J(S)$; hence, $x=0$. This
gives a contradiction. Therefore $S$ is reduced.

Let $a\in R$, there exist $m,n (m>n)$ such that
$\overline{a}^m=\overline{a}^n$ in $S$. Choose $k=n(m-n)$. It is
easy to verify that $p=\overline{a}^{k+1}$ is potent and
$w=\overline{a}-\overline{a}^{k+1}\in N(S)$. Further,
$\overline{a}=p+w=p$ is potent, and so $S$ is potent. Therefore
complete the proof by Lemma 2.7.\hfill$\Box$

\vskip4mm \hspace{-1.8em} {\bf Corollary 2.9.}\ \ {\it Let $R$ be
a ring. Then $R$ is uniquely clean if and only if}
\begin{enumerate}
\item [(1)] {\it $R$ is uniquely $\pi$-clean;}
\vspace{-.5mm}
\item [(2)] {\it $J(R)=\{ x\in R~|~x-1\in U(R)\}.$}
\end{enumerate}
\vspace{-.5mm} {\it Proof.}\ \ Obviously, $J(R)\subseteq \{ x\in
R~|~1-x\in U(R)\}.$ Suppose that $1-x\in U(R)$. Then we have an
idempotent $e\in R$ and an element $u\in J(R)$ such that $x=e+u$
and $ex=xe$ by [11, Theorem 20]. Thus, $1-e=(1-x)+u\in U(R)$, and
so $1-e=1$. This implies that $e=0$; whence, $x=u\in J(R)$.
Therefore $J(R)=\{ x\in R~|~1-x\in U(R)\}$.

Conversely, assume that $(1)$ and $(2)$ hold. In view of Lemma
2.7, $R/J(R)$ is potent. It follows from $J(R)=\{ x\in R~|~x-1\in
U(R)\}$ that $U\big(R/J(R)\big)=\{\overline{1}\}$. Write $p=p^n
(n\geq 2)$ in $R/J(R)$. Then
$\big(1-p^{n-1}+p\big)^{-1}=1-p^{n-1}+p^{n-2}$. Hence,
$p=p^{n-1}$, and so $p^2=p^n=p$. This implies that $R/J(R)$ is
Boolean. Therefore we complete the proof by Lemma 2.1 and [11,
Theorem 20].\hfill$\Box$

\vskip4mm \hspace{-1.8em} {\bf Theorem 2.10.}\ \ {\it Let $R$ be a
ring. Then $R$ is uniquely $\pi$-clean if and only if}
\begin{enumerate}
\item [(1)] {\it For any $a\in R$, there exists an $m\in {\Bbb N}$ and a unique $e\in R$ such that $a^m-e\in J(R)$.}
\vspace{-.5mm}
\item [(2)] {\it $J(R)=\{ x\in R~|~x^m-1\in U(R)~\mbox{for all}~m\in {\Bbb N}\}.$}
\end{enumerate}
\vspace{-.5mm} {\it Proof.}\ \ Suppose that $R$ is uniquely
$\pi$-clean. Let $a\in R$. In view of Theorem 2.8, there exists an
$m\in {\Bbb N}$ and a central idempotent $g\in R$ such that
$a^m-g\in J(R)$. If there exists an idempotent $f\in R$ such that
$a^m-f\in J(R)$, then $g-f=(a^m-f)-(a^m-g)\in J(R)$. Clearly,
$(g-f)^3=g-f$, and so $(g-f)\big(1-(g-f)^2\big)=0$. Thus, $g=f$,
i.e., the uniqueness is verified.

Clearly, $J(R)\subseteq \{ x\in R~|~x^m-1\in U(R)~\mbox{for
all}~m\in {\Bbb N}\}.$ If $x\not\in J(R)$, then $0\neq
xR\nsubseteq J(R)$. In view of Lemma 2.1, $R$ is an exchange ring,
and so there exists an idempotent $0\neq e\in xR$. Write $e=xr$
for a $r\in R$. Choose $a=exe+1-e$. Then we can find some $m\in
{\Bbb N}$ such that $a^m\in R$ is uniquely clean. In addition, $R$
is abelian by Lemma 2.1. Obviously,
$a^m=0+\big(ex^me+1-e\big)=e+\big(e(x^m-1)e+1-e\big)$. If
$x^m-1\in U(R)$, then $0=e$, a contradiction. This implies that
$x^m-1\not\in U(R)$. That is, $x\not\in \{ x\in R~|~x^m-1\in
U(R)~\mbox{for all}~m\in {\Bbb N}\}.$ Therefore $\{ x\in
R~|~x^m-1\in U(R)~\mbox{for all}~m\in {\Bbb N}\}\subseteq J(R)$,
as required.\hfill$\Box$

Conversely, assume that $(1)$ and $(2)$ hold. Let $x\in N(R)$.
Then $x^m-1\in U(R)$ for all $m\in {\Bbb N}$. By hypothesis, we
get $x\in J(R)$. Therefore, every nilpotent element in $R$ is
contained in $J(R)$. Let $e\in R$ be an idempotent, and let $r\in
R$. Then $e+er(1-e)\in R$ is an idempotent. Hence, there exists a
unique $f\in R$ such that $\big(e+er(1-e)\big)-f\in J(R)$. By the
preceding discussion, $\big(e+er(1-e)\big)-e=er(1-e)\in J(R)$. The
uniqueness forces $e=f$. But
$\big(e+er(1-e)\big)-\big(e+er(1-e)\big)\in J(R)$, and so
$e+er(1-e)=f=e$. This shows that $er=ere$. Likewise, $re=ere$.
That is, $er=re$, and then $R$ is abelian. For any $a\in R$, there
exists an $m\in {\Bbb N}$ and a unique $e\in R$ such that
$w:=a^m-e\in J(R)$. Then $a^m=(1-e)+(2e-1+w)$. As $(2e-1)^2=1$, we
see that $2e-1+w\in U(R)$. If there exists an idempotent $f\in R$
such that $a^m-f\in U(R)$, then $e-f=(a^m-f)-(a^m-e)\in U(R)$. One
easily checks that $(e+f-1)(e-f)^2=0$, and therefore $e+f-1=0$.
Thus, $f=1-e$, hence the result.\hfill$\Box$

\vskip4mm \hspace{-1.8em} {\bf Corollary 2.11.}\ \ {\it Let $R$ be
a ring. Then $R$ is uniquely $\pi$-clean if and only if}
\begin{enumerate}
\item [(1)] {\it For any $a\in R$, there exists an $m\in {\Bbb N}$ and a unique $e\in R$ such that $a^m-e\in J(R)$.}
\vspace{-.5mm}
\item [(2)] {\it $N(R)\subseteq J(R)$.}
\end{enumerate}
\vspace{-.5mm} {\it Proof.}\ \ Suppose that $R$ is uniquely
$\pi$-clean. $(1)$ is obvious by Theorem 2.10. Let $a\in N(R)$.
Then $1-a^m\in U(R)$ for all $m\in {\Bbb N}$. It follows by
Theorem 2.10 that $a\in J(R)$. Therefore $N(R)\subseteq J(R)$.

Conversely, assume that $(1)$ and $(2)$ hold. Let $e\in R$, and
let $x\in R$. Then $ex(1-e)\in J(R)$. By hypothesis, we have some
$m\in {\Bbb N}$ such that the expressions
$\big(e+ex(1-e)\big)^m=\big(e+ex(1-e)\big)+0=e+ex(1-e)$ are
unique. This implies that $ex(1-e)=0$, and so $ex=exe$. Likewise,
$xe=exe$. Therefore $R$ is abelian. This yields the result by
Theorem 2.8.\hfill$\Box$

\vskip4mm \hspace{-1.8em} {\bf Corollary 2.12.}\ \ {\it Let $R$ be
a local ring. Then the following are equivalent:}
\begin{enumerate}
\item [(1)] {\it $R$ is uniquely $\pi$-clean.}
\vspace{-.5mm}
\item [(2)] {\it $U(R)=\{ x\in R~|~\exists ~m\in {\Bbb N}~\mbox{such that}~x^m-1\in J(R)\}$.}
\vspace{-.5mm}
\item [(3)] {\it $J(R)=\{ x\in R~|~x^m-1\in U(R)~\mbox{for all}~m\in {\Bbb N}\}.$}
\end{enumerate}
\vspace{-.5mm} {\it Proof.}\ \ $(1)\Rightarrow (3)$ is clear from
Theorem 2.10.

$(3)\Rightarrow (2)$ Obviously, $\{ x\in R~|~\exists ~m\in {\Bbb
N}~\mbox{such that}~x^m-1\in J(R)\}\subseteq U(R)$. For any $x\in
U(R)$, $x\not\in J(R)$. By hypothesis, there exists some $m\in
{\Bbb N}$ such that $x^m-1\not\in U(R)$. As $R$ is local,
$x^m-1\in J(R)$. This implies that $U(R)\subseteq \{ x\in
R~|~\exists ~m\in {\Bbb N}~\mbox{such that}~x^m-1\in J(R)\}$, as
required.

$(2)\Rightarrow (1)$ For any $x\in R$, we see that either $x\in
J(R)$ or $x\in U(R)$. This implies that
$\overline{x}=\overline{0}$ or $\overline{x}^m=\overline{1}$ in
$R/J(R)$. Thus $R/J(R)$ is potent. In light of Lemma 2.7, $R$ is
uniquely $\pi$-clean.\hfill$\Box$

\vskip10mm\bc{\bf 3. FACTORS OF PRIME IDEALS}\ec \vskip4mm The aim
of this section is to characterize uniquely $\pi$-clean rings by
means of prime ideals contains the Jacobson radicals. We use
$J\mbox{-}spec(R)$ to denote the set $\{ P\in
Spec(R)~|~J(R)\subseteq P\}$. Obviously, every maximal ideal is
contained in $J\mbox{-}spec(R)$. Set $$J^*(R)=\bigcap \{
P~|~P~\mbox{is a maximal ideal of}~R\}.$$ We will see that
$J(R)\subseteq J^*(R)$. In general, they are not the same. For
instance, $J(R)=0$ and $J^*(R)=\{ x\in R~|~dim_F(xV)<\infty\}$,
where $R=End_F(V)$ and $V$ is an infinite-dimensional vector space
over a field $F$. Furthermore, we characterize uniquely
$\pi$-clean ring $R$ by means of the radical-like ideal $J^*(R)$.

\vskip4mm \hspace{-1.8em} {\bf Lemma 3.1. [6, Corollary 2.8]}\ \
{\it Let $R$ be a commutative ring. Then the following are
equivalent:}\vspace{-.5mm}
\begin{enumerate}
\item [(1)]{\it $R$ is strongly $\pi$-regular.}
\vspace{-.5mm}
\item [(2)]{\it $R$ is an exchange ring in which every prime ideal of $R$ is maximal.}
\end{enumerate}

\vskip4mm \hspace{-1.8em} {\bf Lemma 3.2.}\ \ {\it Let $R$ be an
abelian exchange ring, and let $x\in R$. Then $RxR=R$ if and only
if $x\in U(R)$.} \vskip2mm\hspace{-1.8em} {\it Proof.}\ \ If $x\in
U(R)$, then $RxR=R$. Conversely, assume that $RxR=R$. As in the
proof of [3, Proposition 17.1.9], there exists an idempotent $e\in
R$ such that $e\in xR$ such that $ReR=R$. This implies that $e=1$.
Write $xy=1$. Then $yx=y(xy)x=(yx)^2$. Hence, $yx=y(yx)x$.
Therefore $1=x(yx)y=xy(yx)xy=yx$, and so $x\in U(R)$. This
completes the proof.\hfill$\Box$

\vskip4mm \hspace{-1.8em} {\bf Theorem 3.3.}\ \ {\it Let $R$ be a
ring. Then $R$ is strongly $\pi$-clean if and only if
}\vspace{-.5mm}
\begin{enumerate}
\item [(1)]{\it $R$ is abelian;}
\vspace{-.5mm}
\item [(2)] {\it Every idempotent lifts modulo $J(R)$;}
\vspace{-.5mm}
\item [(3)]{\it $R/P$ is torsion for all $P\in J\mbox{-}spec(R)$.}
\end{enumerate}\vspace{-.5mm} {\it Proof.}\ \ Suppose $R$ is strongly $\pi$-clean. In view of Lemma 2.1 and
Lemma 2.7, $R$ is an abelian exchange ring, and $R/J(R)$ is
potent. Let $P\in J\mbox{-}spec(R)$. Then $R/J(R)/P/J(R)\cong R/P$
is prime; hence, $P/J(R)$ is a prime ideal of $R/J(R)$. As every
potent ring is commutative, $R/J(R)$ is a commutative
$\pi$-regular ring. It follows from Lemma 3.1 that $P/J(R)$ is a
maximal ideal of $R/J(R)$. We infer that $P$ is a maximal ideal of
$R$.

Clearly, $\overline{R}:=R/P$ is an abelian exchange ring. Since
$P$ is maximal, $R/P$ is simple. For any $0\neq x\in
\overline{R}$, we have $\overline{R}x\overline{R}=\overline{R}$.
By virtue of Lemma 3.2, $x\in U(R/P)$. Hence, $R/P$ is a division
ring. On the other hand, $R/P\cong R/J(R)/P/J(R)$ is potent. Thus,
we have some $m\in {\Bbb N}$ such that $x^{m+1}=x$, and so
$x^m=1$. This implies that $R/P$ is torsion, as required.

Conversely, assume that $(1)-(3)$ hold. Assume that $R$ is not
strongly $\pi$-clean. Set $S=R/J(R)$. In view of Theorem 2.8, $S$
is not periodic. By using Herstein's Theorem, there exists some
$a\in S$ such that $a^m\neq a^{m+1}f(a)$ for any $m\in {\Bbb N}$
and any $f(x)\in {\Bbb Z}[x]$. Let $\Omega=\{ I\lhd
S~|~\overline{a}^{m}\neq \overline{a}^{m+1}f(\overline{a})$ in
$S/I$ for any $m\in {\Bbb N}$ and any $f(x)\in {\Bbb Z}[x]\}$.
Then $\Omega$ is an nonempty inductive. By using Zorn's Lemma,
there exists an ideal $Q$ of $S$ which is maximal in $\Omega$. If
$Q$ is not prime, then there exist two ideals $K$ and $L$ of $R$
such that $K,L\nsubseteq Q$ and $KL\subseteq Q$. By the maximality
of $Q$, we can find some $s,t\in {\Bbb N}$ and some $f(x),g(x)\in
{\Bbb Z}[x]$ such that
$\overline{a}^{s}=\overline{a}^{s+1}f(\overline{a})$ in
$R/\big(K+Q\big)$ and
$\overline{a}^{t}=\overline{a}^{t+1}g(\overline{a})$ in
$R/\big(L+Q\big)$. Thus, $a^{s}-a^{s+1}f(a)\in K+Q$ and
$a^{t}-a^{t+1}g(a)\in L+Q$, and so
$\big(a^{s}-a^{s+1}f(a)\big)\big(a^{t}-a^{t+1}g(a)\big)\in
(K+Q)(L+Q)\subseteq KL+Q\subseteq Q$. In $S/Q$, we have
$\overline{a}^{s+t}=\overline{a}^{s+t+1}h(\overline{a})$ for some
$h(x)\in {\Bbb Z}[x]$. This contradicts the choice of $Q$. Hence,
$Q\in J\mbox{-}spec(R)$. By hypothesis, $R/Q$ is torsion, and so
$R/Q$ is periodic, which is imposable. Therefore $R$ is strongly
$\pi$-clean.\hfill$\Box$

\vskip4mm \hspace{-1.8em} {\bf Corollary 3.4.}\ \ {\it A ring $R$
is uniquely clean if and only if }
\begin{enumerate}
\item [(1)] {\it $R$ is uniquely $\pi$-clean.}
\vspace{-.5mm}
\item [(2)] {\it $R/M\cong {\Bbb Z}_2$ for all maximal ideals
$M$ of $R$.}
\end{enumerate}
\vspace{-.5mm} {\it Proof.}\ \ Suppose $R$ is uniquely clean. Then
$R$ is uniquely $\pi$-clean. $(2)$ is proved by [2, Theorem 2.1].

Conversely, assume that $(1)$ and $(2)$ hold. For all maximal
ideals $M$ of $R$, $1_{R/M}$ is not the sum of two units in $R/M$.
In view of Lemma 2.1, $R$ is an abelian exchange ring, and so it
is clean. Let $x\in R$. Write $x=e_1+u_1=e_1+u_2, e_1=e_1^2,
e_2=e_2^2$ and $u_1,u_2\in U(R)$. If $R\big(1-e_2(1-e_1)\big)R\neq
R$, then there exists a maximal ideal $M$ of $R$ such that
$R\big(1-e_2(1-e_1)\big)R\subseteq M.$ Clearly, $J(R)\subseteq M$.
Hence,
$\overline{x}=\overline{e_1}+\overline{u_1}=\overline{e_2}+\overline{u_2}$
in $R/M$. By Theorem 3.3, $R/M$ is a division ring. This implies
that $\overline{e_i}$ are $\overline{0}$ or $\overline{1}$. If
$\overline{e_1}\neq \overline{e_2}$, then $1_{R/M}$ is the sum of
two units, a contradiction. Therefore we get $e_1-e_2\in M$. This
infers that $e_2(1-e_1)=(e_1-e_2)(e_1-1)\in M$, and so
$1=e_2(1-e_1)+\big(1-e_2(1-e_1)\big)\in M$, a contradiction. As a
result, $R\big(1-e_2(1-e_1)\big)R=R$. As $e_2(1-e_1)\in R$ is an
idempotent, we get $e_2(1-e_1)=0$, and so $e_2=e_2e_1$. Likewise,
$e_1=e_1e_2$. Consequently, $e_1=e_2$, and then $u_1=u_2$.
Therefore $R$ is uniquely clean.\hfill$\Box$

\vskip4mm Let $S(R)$ be the nonempty set of all ideals of a ring
$R$ generated by central idempotents. By Zorn's Lemma, $S(R)$
contains maximal elements. As usual, we say that $R/P$ is a Pierce
stalk if $P$ is a maximal element of the set $S(R)$, and that $P$
is a Pierce ideal. Let $Pier(R)$ be the set of all Pierce ideals
of $R$.

\vskip4mm \hspace{-1.8em} {\bf Proposition 3.5.}\ \ {\it Every
uniquely $\pi$-clean ring is the subdirect product of rings $R_i$,
where each $R_i/J(R_i)$ is torsion.}\vskip2mm\hspace{-1.8em} {\it
Proof.}\ \ Let $R$ be a uniquely $\pi$-clean ring. In view of [9,
Remark 11.2], $\bigcap \{~P~|~P\in Pier(R)\}=0$. Let $\varphi_P:
R\to R/P$ be the natural epimorphism. Then $\bigcap_{P\in
Pier(R)}ker\varphi_P=\bigcap_{P\in Pier(R)}P=0$. Hence, $R$ is the
subdirect product of all $R/P$, where $P\in Pier(R)$. In view of
Lemma 2.1, $R$ is an abelian exchange ring. Let $P\in Pier(R)$.
Then $R/P$ is an exchange ring. As $R$ is indecomposable, we see
that $R/P$ is a local ring. By an argument in [11], $R/P$ is
uniquely $\pi$-clean, and so $R/P/J(R/P)$ is potent from Lemma
2.7, as needed.\hfill$\Box$

\vskip4mm \hspace{-1.8em} {\bf Lemma 3.6.}\ \ {\it Let $R$ be an
abelian exchange ring. Then $J^*(R)=J(R)$.}
\vskip2mm\hspace{-1.8em} {\it Proof.}\ \ Let $M$ be a maximal
ideal of $R$. If $J(R)\nsubseteq M$, then $J(R)+M=R$. Write
$x+y=1$ with $x\in J(R),y\in M$. Then $y=1-x\in U(R)$, an absurd.
Hence, $J(R)\subseteq M$. This implies that $J(R)\subseteq
J^*(R)$. Let $x\in J^*(R)$, and let $r\in R$. If $R(1-xr)R\neq R$,
then we can find a maximal ideal $M$ of $R$ such that
$R(1-xr)R\subseteq M$, and so $1-xr\in M$. It follows that
$1=xr+(1-xr)\in M$, which is imposable. Therefore $R(1-xr)R=R$. In
light of Lemma 3.2, $1-xr\in U(R)$, and then $x\in J(R)$. This
completes the proof.\hfill$\Box$

\vskip4mm \hspace{-1.8em} {\bf Theorem 3.7.}\ \ {\it Let $R$ be a
ring. Then $R$ is uniquely $\pi$-clean if and only if}
\begin{enumerate}
\item [(1)] {\it $R$ is an exchange ring;}
\vspace{-.5mm}
\item [(2)] {\it $R/J^*(R)$ is potent and every idempotent uniquely lifts modulo $J^*(R)$.}
\end{enumerate}
\vspace{-.5mm} {\it Proof.}\ \ Suppose $R$ is uniquely
$\pi$-clean. Then $R$ is an abelian exchange ring by Lemma 2.1. In
view of Lemma 3.6, $J^*(R)=J(R)$. It follows from Lemma 2.7 that
$R/J^*(R)$ is potent. Let $e-e^2\in J(R)$. Then we can find an
idempotent $f\in R$ such that $e-f\in J(R)$. Since
$(e-f)^2\big(1-(e-f)\big)=0$, we get $e=f$, as desired.

Conversely, assume that $(1)$ and $(2)$ hold. Let $e\in R$ be an
idempotent, and let $r\in R$. Then $\overline{er(1-e)}\in
R/J^*(R)$ is potent. This implies that
$\overline{er(1-e)}=\overline{0}$, and so $er(1-e)\in J^*(R)$.
Since $e-e, e-\big(e+er(1-e)\big)\in J^*(R)$, by the uniqueness,
we get $e=e+er(1-e)$, and so $er=ere$. Likewise, $re=ere$; hence
that $er=re$. Thus, $R$ is abelian. In light of Lemma 3.6,
$J^*(R)=J(R)$. Therefore we complete the proof, in terms of Lemma
2.7.\hfill$\Box$

\vskip4mm \hspace{-1.8em} {\bf Corollary 3.8.}\ \ {\it Let $R$ be
a ring which have finitely many maximal ideals. Then $R$ is
uniquely $\pi$-clean if and only if}
\begin{enumerate}
\item [(1)] {\it $R$ is an exchange ring;}
\vspace{-.5mm}
\item [(2)] {\it $R/J^*(R)$ is the direct sum of finitely many torsion rings and every idempotent uniquely lifts modulo $J^*(R)$.}
\end{enumerate}
\vspace{-.5mm} {\it Proof.}\ \ $(1)\Rightarrow (2)$ Let $M$ be a
maximal ideal of $R$. As in the proof of Lemma 3.6, we see that
$J(R)\subseteq M$. This shows that $M\in J\mbox{-}spec(R)$.
Therefore $R/M$ is torsion by Theorem 3.3. Since $R$ has finitely
many maximal ideals $M_1,\cdots ,M_n$, we see that $R/J^*(R)\cong
R/M_1\oplus \cdots \oplus R/M_n$. Therefore $R/J^*(R)$ is the
direct sum of finitely many torsion rings.

$(2)\Rightarrow (1)$ As every torsion ring is potent, we see that
$R/J^*(R)$ is potent. Therefore we complete the proof, by Theorem
3.7.\hfill$\Box$

\vskip4mm \hspace{-1.8em} {\bf Theorem 3.9.}\ \ {\it Let $R$ be a
ring. Then $R$ is uniquely $\pi$-clean if and only if}
\begin{enumerate}
\item [(1)] {\it For any $a\in R$, there exists an $m\in {\Bbb N}$ and a unique $e\in R$ such that $a^m-e\in J^*(R)$.}
\vspace{-.5mm}
\item [(2)] {\it $J^*(R)=\{ x\in R~|~x^m-1\in U(R)~\mbox{for all}~m\in {\Bbb N}\}.$}
\end{enumerate}
\vspace{-.5mm} {\it Proof.}\ \ One direction is obvious by Lemma
3.6 and Theorem 2.10.

Conversely, assume that $(1)$ and $(2)$ hold. Let $x\in N(R)$.
Then $x^m-1\in U(R)$ for all $m\in {\Bbb N}$. By hypothesis, we
have $x\in J^*(R)$, and so $N(R)\subseteq J^*(R)$. Let $e\in R$ be
an idempotent, and let $r\in R$. Then $e+er(1-e)+0=e+er(1-e)$ with
$0,er(1-e)\in J^*(R)$. By the uniqueness, we get $er=ere$.
Similarly, we have $re=ere$. That is, $er=re$. We infer that $R$
is abelian. For any $a\in R$, there exists an $m\in {\Bbb N}$ and
a unique $e\in R$ such that $w:=a^m-e\in J^*(R)$. Then
$a^m=(1-e)+(2e-1+w)$. But $2e-1+w=(1-2e)\big((1-2e)w-1\big)\in
U(R)$, by $(2)$. If there exists an idempotent $f\in R$ such that
$a^m-f\in U(R)$, then
$e-f=(a^m-f)-(a^m-e)=(a^m-f)\big(1-(a^m-f)^{-1}(a^m-e)\big)\in
U(R)$. It follows from $(e+f-1)(e-f)^2=0$ that $f=1-e$, and we are
through.\hfill$\Box$

\vskip4mm Let $P(R)$ be the intersection of all prime ideals of
$R$, i.e., $P(R)$ is the prime radical of R. As is well known,
$P(R)$ is the intersection of all minimal prime ideals of $R$.

\vskip4mm \hspace{-1.8em} {\bf Corollary 3.10.}\ \ {\it Let $R$ be
a uniquely $\pi$-clean in which every prime ideal is maximal. Then
$$P(R)=\{ x\in R~|~x^m-1\in U(R)~\mbox{for all}~m\in {\Bbb N}\}.$$}
\hspace{-1.0em} {\it Proof.}\ \ As every maximal ideal is prime,
we deduce that $J^*(R)=P(R)$, and therefore we complete the proof
by Theorem 3.9.\hfill$\Box$

\vskip10mm\bc{\bf 4. CERTAIN CLASSES}\ec \vskip4mm In this section
we investigate certain classes of uniquely $\pi$-clean rings. So
as to construct more examples of uniquely clean rings, we recall
the concept of ideal-extensions. Let $R$ be a ring with an
identity and $S$ be a ring (not necessary unitary), and let $S$ be
a $R$-$R$-bimodule in which $(s_1s_2)r=s_1(s_2r),
r(s_1s_2)=(rs_1)s_2$ and $(s_1r)s_2=s_1(rs_2)$ for all $s_1,s_2\in
S, r\in R$. The ideal extension $I(R;S)$ is defined to be the
additive abelian group $R\oplus S$ with multiplication $
(r_1,s_1)(r_2,s_2)=(r_1r_2,s_1s_2+r_1s_2+s_1r_2)$ (see [4] and
[10]). We start this section by examining when an ideal extension
is uniquely $\pi$-clean.

\vskip4mm \hspace{-1.8em} {\bf Theorem 4.1.}\ \ {\it The
ideal-extension $I(R;S)$ is uniquely $\pi$-clean and $S$ is
idempotent-free if and only if }\vspace{-.5mm}
\begin{enumerate}
\item [(1)]{\it $R$ is uniquely $\pi$-clean;}
\vspace{-.5mm}
\item [(2)]{\it If $e=e^2\in R$, then $es=se$ for all $s\in S$;}
\vspace{-.5mm}
\item [(3)]{\it If $s\in S$, then there exists some $s'\in S$ such that $ss'=s's$ and $s+s'+ss'=0$.}
\end{enumerate}\vspace{-.5mm}  {\it Proof.}\ \ Assume that $(1)-(3)$
hold. Let $e\in S$ be an idempotent. Then $(-e)+s'+(-e)s'=0$ for
some $s'\in S$. Hence, $(1-e)(1+s')=1$, and so $e=0$. That is, $S$
is idempotent-free. Let $(a,s)\in I(R;S)$. Then $a\in R$ is
uniquely $\pi$-clean. Thus, we have some $n\in {\Bbb N}$ such that
$a^n\in R$ is uniquely clean. Write $a^n=e+u,e=e^2\in R, u\in
U(R)$. Hence, $(a,s)^n=(a^n,x)=(e,0)+(u,x)$ for some $x\in S$.
Clearly, $(e,0)^2=(e,0)$ and $(u,x)^{-1}=(u^{-1},zu^{-1})$ for a
$z\in S$. Write $(a,s)^n=(f,y)+(v,w), (f,y)^2=(f,y)$ and $(v,w)\in
U\big(I(R;S)\big)$. Then $f=f^2\in R,y=0$ and $v\in U(R)$.
Further, $a^n=f+v$ and $x=w$. This implies that $f=e,v=u$, and so
$(f,y)=(e,0),(v,w)=(u,x)$. As a result, $(a,s)\in I(R;S)$ is
uniquely $\pi$-clean, and so $I(R;S)$ is uniquely $\pi$-clean.

Assume that $I(R;S)$ is uniquely $\pi$-clean and $S$ is
idempotent-free. Then $R$ is uniquely $\pi$-clean. Let $e=e^2\in
R$ and $s\in S$. In view of Lemma 2.1, $(e,0)=(e,0)^2\in I(R;S)$
is central. Hence, $(e,0)(0,s)=(0,s)(e,0)$, and so $es=se$. For
any $s\in S$, there exists some $n\in {\Bbb N}$ such that
$(1,s)^n\in I(R;S)$ is uniquely clean. Write
$(1,s)^n=(1,x)=(f,y)+(u,v)$ where $x\in S, (f,y)\in I(R;S)$ is an
idempotent and $(u,v)\in I(R;S)$ is a unit. Clearly, $f=0$, and so
$y=0$. This implies that $x=y+v=v$; hence, $(1,x)\in I(R;S)$ is a
unit. Further, $(1,s)\in I(R;S)$ is a unit. Write
$(1,s)^{-1}=(1,s')$ for a $s'\in S$. Then $ss'=s's$ and
$s+s'+ss'=0$, hence the result.\hfill$\Box$

\vskip4mm \hspace{-1.8em} {\bf Corollary 4.2.}\ \ {\it Let $R$ be
uniquely $\pi$-clean. Then $S=\{ (a_{ij})\in
T_n(R)~|~a_{11}=\cdots =a_{nn}\}$ is uniquely $\pi$-clean.}
\vskip2mm\hspace{-1.8em} {\it Proof.}\ \ Let $T=\{ (a_{ij})\in
T_n(R)~|~a_{11}=\cdots =a_{nn}=0\}$. Then $S\cong I(R;T)$. Then
the result follows by Theorem 4.1.\hfill$\Box$

\vskip4mm A ring $R$ is called potently $J$-clean if for any $a\in
R$ there exists a potent $p\in R$ such that $a-p\in J(R)$. We
shall show that such rings form a subclass of uniquely $\pi$-clean
rings.

\vskip4mm \hspace{-1.8em} {\bf Lemma 4.3.}\ \ {\it Every potently
$J$-clean ring is an exchange ring.} \vskip2mm\hspace{-1.8em} {\it
Proof.}\ \ Let $R$ be a potently $J$-clean ring. Then $R/J(R)$ is
potent, and so it is an exchange ring. Let $\overline{e}\in
R/J(R)$ be an idempotent. Then we have a potent $p\in R$ such that
$w:=e-p\in J(R)$. Write $p=p^n$ for some $n\geq 2$. Then
$p^{n-1}\in R$ is an idempotent. Moreover, $e=p+w$, and so
$e^{n-1}=p^{n-1}+v$ for some $v\in J(R)$. But $e-e^{n-1}\in J(R)$.
Hence, $e-p^{n-1}=\big(e-e^{n-1}\big)+\big(e^{n-1}-p^{n-1}\big)\in
J(R)$. So idempotents can be lifted modulo $J(R)$. In light of [9,
Theorem 29.2], $R$ is an exchange ring.\hfill$\Box$

\vskip4mm \hspace{-1.8em} {\bf Theorem 4.4.}\ \ {\it Every abelian
potently $J$-clean ring is uniquely $\pi$-clean.}
\vskip2mm\hspace{-1.8em} {\it Proof.}\ \ Let $R$ be a abelian
potently $J$-clean ring. Then $R$ is an exchange ring by Lemma
4.3. Thus, every idempotent in $R$ lifts modulo $J(R)$. For any
$a\in R$, there exists a potent $p\in R$ such that $a-p\in J(R)$.
This implies that $\overline{a}\in R/J(R)$ is potent, and so
$R/J(R)$ is potent. According to Lemma 2.7, $R$ is uniquely
$\pi$-clean.\hfill$\Box$

\vskip4mm \hspace{-1.8em} {\bf Corollary 4.5.}\ \ {\it Let $R$ be
abelian. If for any sequence of elements $\{ a_i\}\subseteq R$
there exists a $k\in {\Bbb N}$ and $n_1,\cdots ,n_k\geq 2$ such
that $(a_1-a_1^{n_1})\cdots (a_k-a_k^{n_k})=0$, then $R$ is
uniquely $\pi$-clean.} \vskip2mm\hspace{-1.8em} {\it Proof.}\ \
For any $a\in R$, we have a $k\in {\Bbb N}$ and $n_1,\cdots
,n_k\geq 2$ such that $(a-a^{n_1})\cdots (a-a^{n_k})=0$. This
implies that $a^k=a^{k+1}f(a)$ for some $f(t)\in {\Bbb Z}[t]$. By
Herstein's Theorem, $R$ is periodic. Therefore every element in
$R$ is the sum of a potent element and a nilpotent element.

Clearly, $R/J(R)$ is isomorphic to a subdirect product of some
primitive rings $R_i$. Case 1. There exists a subring $S_i$ of
$R_i$ which admits an epimorphism $\phi_i: S_i\to M_{2}(D_i)$
where $D_i$ is a division ring. Case 2. $R_i\cong M_{m_i}(D_i)$
for a division ring $D_i$. Clearly, the hypothesis is inherited by
all subrings, all homomorphic images and all corners of $R$, we
claim that, for any sequence of elements $\{ a_i\}\subseteq
M_2(D_i)$ there exists $s\in {\Bbb N}$ and $m_1,\cdots ,m_s\geq 2$
such that $(a_1-a_1^{m_1})\cdots (a_s-a_s^{m_s})=0$. Choose
$a_i=e_{12}$ if $i$ is odd and $a_i=e_{21}$ if $i$ is even. Then
$(a_1-a_1^{m_1})(a_2-a_2^{m_2})\cdots (a_s-a_s^{m_s})=a_1a_2\cdots
a_s\neq 0$, a contradiction. This forces $m_i=1$ for all $i$. We
infer that all $R_i$ is reduced, and then so is $R/J(R)$. If $a\in
N(R)$, we have some $n\in {\Bbb N}$ such that $a^n=0$, and thus
$\overline{a}^n=0$ is $R/J(R)$. Hence, $\overline{a}\in
J\big(R/J(R)\big)=0$. This implies that $a\in J(R)$, and so
$N(R)\subseteq J(R)$. Therefore $R$ is potently $J$-clean, hence
the result by Theorem 4.4.\hfill$\Box$

\vskip4mm Recall that an element $a\in R$ is uniquely nil clean
provided that there exists a unique idempotent $e\in R$ such that
$a-e\in N(R)$ [2].

\vskip4mm \hspace{-1.8em} {\bf Lemma 4.6.}\ \ {\it Let $R$ be a
ring. Then the following are equivalent:}
\begin{enumerate}
\item [(1)] {\it $R$ is an abelian periodic ring.}
\vspace{-.5mm}
\item [(2)] {\it For any $a\in R$, there exists some $m\in {\Bbb N}$ such
that $a^m\in R$ is uniquely nil clean.}
\end{enumerate}
\vspace{-.5mm}  {\it Proof.}\ \ $(1)\Rightarrow (2)$ Let $a\in R$.
Since $R$ is periodic, there exists some $m\in {\Bbb N}$ such that
$a^m\in R$ is an idempotent. Write $a^m=e+w$ where $e=e^2\in R$
and $w\in N(R)$. Then $a^m-e=w\in N(R)$. As $R$ is abelian, we see
that $\big(a^m-e\big)^3=a^m-e$. Thus,
$\big(a^m-e\big)\big(1-(a^m-e)^2\big)=1$, and so $a^m=e$, as
required.

$(2)\Rightarrow (1)$ Let $e\in R$ be an idempotent and $r\in R$.
Choose $a=e+er(1-e)$. Then we can find some $m\in {\Bbb N}$ such
that $a^m\in R$ is uniquely nil clean. As
$a=a^m=e+er(1-e)=\big(e+er(1-e)\big)+0$, by the uniqueness, we get
$er(1-e)=0$, and so $er=ere$. Likewise, $re=ere$, and so $er=re$.
Therefore $R$ is abelian. Let $a\in R$. Then there exists some
$n\in {\Bbb N}$ such that $a^n=f+u$, where $f=f^2\in R$ and $u\in
N(R)$. Hence, $a^{2n}=f+v$ for a $v\in N(R)$ and $uv=vu$. This
shows that $a^n-a^{2n}=u-v\in N(R)$. Thus, we have a $k\in {\Bbb
N}$ such that $a^{nk}=a^{nk+1}f(a)$ for some $f(t)\in {\Bbb
Z}[t]$. In light of Herstein's Theorem, $R$ is
periodic.\hfill$\Box$

\vskip4mm \hspace{-1.8em} {\bf Theorem 4.7.}\ \ {\it Let $R$ be a
ring. Then the following are equivalent:}\vspace{-.5mm}
\begin{enumerate}
\item [(1)]{\it $R$ is uniquely $\pi$-clean and $J(R)$ is nil.}
\vspace{-.5mm}
\item [(2)] {\it $R$ is an abelian periodic ring.}
\vspace{-.5mm}
\item [(3)]{\it For any $a\in R$, there exists some $m\in {\Bbb N}$ and a unique idempotent $e\in R$ such that $a^m-e\in P(R)$.}
\end{enumerate}\vspace{-.5mm} {\it Proof.}\ \ $(1)\Rightarrow (2)$
In view of Lemma 2.7, $R/J(R)$ is potent. Let $a\in R$. Then
$\overline{a}=\overline{a^m} (m\geq 2)$, and so
$\overline{a}^{m-1}$ is an idempotent. As $J(R)$ is nil, every
idempotent lifts modulo $J(R)$. Hence, we can find an idempotent
$e\in R$ such that $a^{m-1}=e+w$, where $w\in J(R)$ is nilpotent.
Write $a^{m-1}=f+v$ with $f=f^2\in R$ and $v\in N(R)$. In view of
Lemma 2.1, $R$ is abelian. Then $e-f=v-w\in N(R)$, as $vw=wv$. It
follows from $e-f=(e-f)^3$ that $e=f$. Thus, proving $(2)$ by
Lemma 4.6.

$(2)\Rightarrow (3)$ For any $a\in R$, by Lemma 4.6, there exists
some $m\in {\Bbb N}$ such that $a^m$ is uniquely nil clean. Write
$a^m=e+w$ with $e=e^2$ and $w\in N(R)$. In view of [1, Theorem 2],
$N(R)$ forms an ideal of $R$. Therefore $N(R)=P(R)$, as required.

$(3)\Rightarrow (1)$ Let $e\in R$ be an idempotent, and let $r\in
R$. Then we have an idempotent $f\in R$ such that $er(1-e)=f+w$
for a $w\in P(R)$. Hence,
$1-f=1-er(1-e)+w=\big(1-er(1-e)\big)\big(1+(1+er(1-e))w\big)\in
U(R)$. We infer that $f=0$, and so $er(1-e)=w\in P(R)$. But we
have a unique expression $e+er(1-e)=e+er(1-e)+0$ where $er(1-e),
0\in P(R)$. By the uniqueness, we get $e=e+er(1-e)$, and so
$er=ere$. Similarly, $re=ere$. Therefore $er=re$, i.e., $R$ is
abelian.

Let $x\in J(R)$. Write $x=h+v$ with $h=h^2\in R, v\in P(R)$. Then
$h=x-v\in J(R)$; hence that $h=0$. It follows that $J(R)=P(R)$.
Accordingly, for any $a\in R$, there exists some $m\in {\Bbb N}$
and a unique idempotent $e\in R$ such that $a^m-e\in J(R)$.

If $x\in N(R)$, then we have an idempotent $g\in R$ and a $u\in
P(R)$ such that $x=g+u$, and so $g=x-u$. As $R$ is abelian, we see
that $xu=ux$, and then $g\in N(R)$. This shows that $g=0$.
Consequently, $x=u\in P(R)\subseteq J(R)$. We infer that
$N(R)\subseteq J(R)$. In light of Corollary 2.11, we complete the
proof.\hfill$\Box$

\vskip4mm As every finite ring is periodic, it follows from
Theorem 4.7 that every finite commutative ring is uniquely
$\pi$-clean, e.g., ${\Bbb Z}_{n}[\alpha]=\{ a+b\alpha~|~a,b\in
{\Bbb Z}_{n}, \alpha=-\frac{1}{2}+\frac{\sqrt{3}}{2}i, i^2=-1\}$.

\vskip4mm \hspace{-1.8em} {\bf Corollary 4.8.}\ \ {\it Let $R$ be
a ring. Then the following are equivalent:}\vspace{-.5mm}
\begin{enumerate}
\item [(1)]{\it $R$ is uniquely $\pi$-clean and $J(R)$ is nil.}
\vspace{-.5mm}
\item [(2)] {\it For any $a\in R$, there exists some $m\in {\Bbb N}$ and a central idempotent $e\in R$ such that $a^m-e\in P(R)$.}
\end{enumerate}\vspace{-.5mm} {\it Proof.}\ \ $(1)\Rightarrow (2)$ This is obvious, in
terms of Theorem 4.7 and Lemma 2.1.

$(2)\Rightarrow (1)$ For any $a\in R$, there exists some $m\in
{\Bbb N}$ and a central idempotent $e\in R$ such that $a^m-e\in
P(R)$. Write $a^m-f\in P(R)$ for an idempotent. Then
$e-f=(a^m-f)-(a^m-e)\in P(R)$. As $(e-f)^3=e-f$, we conclude that
$e=f$, and we are through by Theorem 4.7. \hfill$\Box$

\vskip4mm Let $n\geq 2$ be a fixed integer. Following Yaqub, a
ring $R$ is said to be generalized $n$-like provided that for any
$a,b\in R$, $(ab)^n-ab^n-a^nb+ab=0$ ([7-8]).

\vskip4mm \hspace{-1.8em} {\bf Corollary 4.9.}\ \ {\it Every
generalized $n$-like ring is uniquely $\pi$-clean.}
\vskip2mm\hspace{-2.0em} {\it Proof.}\ \ Let $a\in R$. Then
$a^{2n}-2a^{n+1}+a^2=0$, and so $(a-a^n)^2=0$. Thus, $a-a^n\in
N(R)$. Hence, $a^m=a^{m+1}f(a)$ for some $f(t)\in {\Bbb Z}[t]$.
Accordingly, $R$ is periodic by Herstein's Theorem.

Let $e,f\in R$. Then there exist some $m,n\geq 2$ such that
$$\begin{array}{c}
\big((1-e)f\big)^me=\big((1-e)fe\big)^m-(1-e)fe+(1-e)fe=0;\\
\big((1-e)f\big)^n=(1-e)f+(1-e)f-(1-e)f=(1-e)f.
\end{array}
$$ Reiterating in the last, we get
$(1-e)f=\big((1-e)f\big)^{n+m}$, and so $(1-e)fe=0$. Hence,
$fe=efe$. Likewise, $ef=efe$. Therefore $ef=fe$. We infer that $R$
is abelian.

Therefore we conclude that $R$ is uniquely $\pi$-clean, in terms
of Theorem 4.7.\hfill$\Box$

\vskip4mm Let $R=\{
\left( \begin{array}{ccc} x&y&z\\
0&x^2&0\\
0&0&x
\end{array}
\right)~|~x,y,z\in GF(4)\}$. It is easy to check that for each
$a\in R$, $a^7=a$ or $a^7=a^2=0$. Therefore $R$ is a generalized
$7$-like ring. By Corollary 4.9, $R$ is uniquely $\pi$-clean which
is a noncommutative periodic ring.

\vskip20mm \bc{\bf REFERENCES}\ec \vskip4mm {\small \re{1} A.
Badawi, On abelian $\pi$--regular rings, {\it{Comm. Algebra,}}
{\bf{25}}(1997), 1009--1021.

\re{2} H. Chen, On uniquely clean rings, {\it Comm. Algebra}, {\bf
39}(2011), 189--198.

\re{3} H. Chen, {\it Rings Related Stable Range Conditions},
Series in Algebra 11, World Scientific, Hackensack, NJ, 2011.

\re{4} V.A. Hiremath and S. Hegde, Using ideals to provide a
unified approach to uniquely clean rings, {\it J. Aust. Math.
Soc.}, {\bf 96}(2014), 258--274.

\re{5} M.T. Kosan; T.K. Lee and Y. Zhou, When is every matrix over
a division ring a sum of an idempotent and a nilpotent? {\it
Linear Algebra Appl.}, {\bf 450}(2014), 7--12.

\re{6} D. Lu and W. Yu, On prime spectrums of commutative rings,
{\it Comm. Algebra}, {\bf 34}(2006), 2667--2672.

\re{7} H.G. Moore, Generalized $n$-like rings and commutativity,
{\it Canad. Math. Bull.}, {\bf 23} (1980), 449--452.

\re{8} H. Tominaga and A. Yaqub, On generated $n$-like rings and
related rings, {\it Math. J. Okayama Univ.}, {\bf 23}(1981),
199--202.

\re{9} A.A. Tuganbaev, {\it Rings Close to Regular}, Kluwer
Academic Publishers, Dordrecht, Boston, London, 2002.

\re{10} W.K. Nicholson and Y. Zhou, Rings in which elements are
uniquely the sum of an idempotent and a unit, {\it Glasgow Math.
J.}, {\bf 46}(2004), 227--236.

\re{11} Y. Zhou, Some recent work on clean rings, preprint, 2011.
\end{document}